# Nonlinear Split Ordered Variational Inequality Problems

Jinlu Li
Department of Mathematics
Shawnee State University
Portsmouth, Ohio 45662
USA

## Abstract

The concept of nonlinear split ordered variational inequality problems on partially ordered vector spaces is a natural extension of linear split vector variational inequality problems on Banach spaces. The results about nonlinear split ordered variational inequality problems are immediately applied to solving nonlinear split vector optimization problems. In this paper, we prove the solvability of some nonlinear split ordered variational inequality problems by using some fixed point theorems on partially ordered spaces, in which the considered mappings may not be required to have any type of continuity and they just satisfy some order-monotonic conditions. As applications of the results about nonlinear split ordered variational inequality problems, we study solvability of some nonlinear split variational inequality problems and linear split variational inequality problems on partially ordered vector spaces and partially ordered Banach spaces.

## 1. Introduction

Let $X$ and $Y$ be Hilbert spaces with dual spaces $X^*$ and $Y^*$, respectively. For any given nonempty, closed and convex subsets $C \subseteq X$ and $D \subseteq Y$, given operators $f: C \to X^*$, $g: D \to Y^*$ and a bounded linear operator $A: X \to Y$, the linear (in view of the linearity of $A$) split variational inequality problem (or simply called the split variational inequality problem) associated with $f$, $g$, $A$, $C$, $D$, denoted by SVI($f$, $C$, $A$, $g$, $D$), is formulated as follows (see [4–7]):

$$\text{find } x_* \in C \text{ such that } \langle f(x_*), x - x_* \rangle \geq 0 \text{ for all } x \in C, \tag{1}$$

and such that

$$y_* = Ax_* \in D \text{ solves } \langle g(y_*), y - y_* \rangle \geq 0 \text{ for all } y \in D. \tag{2}$$

This concept has been extended to Banach spaces (see [14–15], [19–20]). It is clear that split variational inequality problems naturally generalized variational inequality problems on the corresponding underlying spaces.

According to my known knowledge, in the publications about split variational inequality problems, all authors provide nice iterated sequences to approximate the solutions for some considered problems under a very strong condition: Assume that the considered split variational inequality problem has a solution (see [2–10], [17], [19–22]). In these publications, no author studied the existence of solutions to any considered split variational inequality problem on

Hilbert spaces or on Banach spaces, in which the assumption of the existence of solutions has been repeatedly applied in their proofs for the convergence of the algorithms provided.

In [14], an existence theorem for solutions to some split variational inequality problems is proved by using the Fan–KKM Theorem (see [11], [18]), which is extended in [15] to the solvability of split vector variational inequality problems. In order to prove these results, new concepts are introduced which are named as convexity direction preserving property for the considered mappings.

Similar to the applications of split variational inequality problems to optimization problems (see [1], [13], [16]), split (vector) variational inequality problems have been applied to solving some split (vector) optimization problems ([4], [9], [14–15]).

In this paper, we generalize (linear) split vector variational inequality problems to nonlinear split ordered variational inequality problems, in which the operator $A$ appeared in SVI($f$, $C$, $A$, $g$, $D$) will be extended from linear to nonlinear. It is noticeable that the key point for the solvability of split (vector) variational inequality problems is that, in addition to the conditions for the operators' $f$ and $g$, some conditions for the operator $A$ must be investigated. One of the main goals of this paper is to find some applicable conditions for $A$, $f$, and $g$, which are different from the convexity direction preserving property used in [14–15], to assure the existence of solutions to some nonlinear split ordered variational inequality problems. Consequently, the problem SVI($f$, $C$, $A$, $g$, $D$) defined in (1) and (2) and some split vector variational inequality problems can be solved.

## 2. Preliminaries

In this section, we review some concepts about partially ordered sets and recall some fixed point theorem on posets. For details, the readers are referred to [12–13].

A vector space $X$ equipped with a partially order $\succcurlyeq$ on it is called a partially ordered vector space if the following order linearity properties hold:

1. $x \succcurlyeq y$ implies $x + z \succcurlyeq y + z$, for all $x, y, z \in X$;
2. $x \succcurlyeq y$ implies $\alpha x \succcurlyeq \alpha y$, for all $x, y \in X$ and $\alpha \geq 0$.

A topological space $(X, \tau)$ equipped with a partially order $\succcurlyeq$ on it is called a partially ordered topological space, denoted by $(X, \tau, \succcurlyeq)$, if

3. The topology $\tau$ is natural with respect to the partial order $\succcurlyeq$, that is, for every $u \in B$, all the following $\succcurlyeq$-intervals are $\tau$-closed
$$[u) = \{x \in B: x \succcurlyeq u\} \text{ and } (u] = \{x \in B: x \preccurlyeq u\}. \tag{3}$$

A topological vector space $(X, \tau)$ equipped with a partially order $\succcurlyeq$ on it is called a partially ordered topological vector space, denoted by $(X, \tau, \succcurlyeq)$, if it is both of partially ordered vector space and partially ordered topological space, that is, all above three conditions are satisfied.

In particular, let $(X, \|\cdot\|, \succcurlyeq)$ be a partially ordered Banach space with the weak topology $\omega$. Since, for every $u \in B$, all the $\succcurlyeq$-intervals in (3) are convex, then these $\succcurlyeq$-intervals are $\|\cdot\|$-closed if and only if, they are weakly $\omega$-closed. It implies that, the norm $\|\cdot\|$-topology is natural with respect to the partial order $\succcurlyeq$, if and only if, the weak topology $\omega$ is natural with respect to the partial order $\succcurlyeq$. That is,

$$(X, \|\cdot\|, \succcurlyeq) \text{ is a partially ordered Banach space}$$

if, and only if

$$(X, \omega, \succcurlyeq) \text{ is a partially ordered topological space.} \tag{4}$$

A nonempty subset $A$ of a poset $(X, \succcurlyeq)$ is said to be chain-complete, if every chain contained in $A$ admits its smallest $\succcurlyeq$-upper bound in $A$. $A$ is said to be universally inductive in $X$ if, for any given chain $\{x_\alpha\} \subseteq X$ satisfying that every element $x_\beta \in \{x_\alpha\}$ has an upper $\succcurlyeq$-cover in $A$, then the chain $\{x_\alpha\}$ has an upper $\succcurlyeq$-bound in $A$.

Examples of chain-complete posets and universally inductive posets:

1. Every non-empty compact subset of a partially ordered Hausdorff topological space is both chain-complete and universally inductive.
2. Every non-empty bounded closed and convex subset of a partially ordered reflexive Banach space is chain-complete and universally inductive.
3. Every inductive subset of a chain-complete poset with a finite number of maximal elements is universally inductive.
4. Every nonempty closed inductive subset of a regular partially ordered Banach space is chain-complete.

We recall a fixed point theorem from [13] that will be used in the following contents.

**Fixed Point Theorem A** (Theorem 3.2 in [13]). *Let $(P, \succcurlyeq)$ be a chain complete poset and let $F$: $P \to 2^P \setminus \{\emptyset\}$ be a set-valued mapping satisfying the following three conditions*:

   A1. *$F$ is $\succcurlyeq$-increasing upward*;
   A2. *$(F(x), \succcurlyeq)$ is universally inductive, for every $x \in P$*;
   A3. *There is an element $y_*$ in $P$ and $v_* \in F(y_*)$ with $y_* \preccurlyeq v_*$.*

*Then*

   (i)    *$(\mathcal{F}(F), \succcurlyeq)$ is a nonempty inductive poset*;
   (ii)   *$(\mathcal{F}(F) \cap [y_*), \succcurlyeq)$ is a nonempty inductive poset; and $F$ has an $\succcurlyeq$-maximal fixed point $x^*$ with $x^* \succcurlyeq y_*$.*

**Fixed Point Theorem B.** *Let $(X, \tau, \succcurlyeq)$ be a partially ordered compact Hausdorff topological space. Let $F$: $X \to 2^X \setminus \{\emptyset\}$ be a set-valued mapping satisfying the following three conditions*:

   A1. *$F$ is $\succcurlyeq$-increasing upward*;
   A2. *$F(x)$ is a nonempty closed subset, for every $x \in X$*;
   A3. *There is an element $y_*$ in $X$ and $v_* \in F(y_*)$ with $y_* \preccurlyeq v_*$.*
*Then*

(i) $(\mathcal{F}(F), \succcurlyeq)$ *is a nonempty inductive poset*;

  (ii) $(\mathcal{F}(F) \cap [y_*), \succcurlyeq)$ *is a nonempty inductive poset, and F has an $\succcurlyeq$-maximal fixed point $x^*$ with $x^* \succcurlyeq y_*$.*

## 3. Nonlinear (*A* is nonlinear) split ordered variational inequality problems on vector spaces

3.1 Definitions of Nonlinear split ordered variational inequality problems on vector spaces

For any vector spaces $X$ and $Y$, let $L(X, Y)$ denote the space of linear operators from $X$ to $Y$. If $X$ and $Y$ are topological vector spaces, then let $\mathcal{L}(X, Y)$ denote the space of linear continuous operators from $X$ to $Y$.

**Definition 3**. Let $X$ and $Y$ be vector spaces. Let $(U, \succcurlyeq^U)$ and $(V, \succcurlyeq^V)$ be partially ordered vector spaces. For any given nonempty subsets $C \subseteq X$ and $D \subseteq Y$, given operators $f: C \to L(X, U)$, $g: D \to L(Y, V)$ and an nonlinear operator $A: X \to Y$, the nonlinear (based on the nonlinearity of the operator $A$) split ordered variational inequality problem associated with $f$, $g$, $A$, $C$, $D$, $U$, and $V$ (denoted by NSOVI($f, C, U, A, g, D, V$)) is formulated as follows:

$$\text{find } x_* \in C \text{ such that } f(x_*)(x - x_*) \succcurlyeq^U 0, \text{ for all } x \in C, \tag{5}$$

and such that

$$y_* = Ax_* \in D \text{ solves } g(y_*)(y - y_*) \succcurlyeq^V 0, \text{ for all } y \in D. \tag{6}$$

This point $x^* \in C$, (or the pair $(x_*, y_*) = (x_*, Ax_*) \in C \times D$) is called a solution to this problem NSOVI($f, C, U, A, g, D, V$) and the solution set is denoted by $\mathcal{S}(f, A, g)$.

**Remarks**. One can consider the nonlinear split ordered variational inequality problem NSOVI($f, C, U, A, g, D, V$) to be formulated as follows:

$$\text{find } x_* \in C \text{ such that } f(x_*)(x - x_*) \not\prec^U 0, \text{ for all } x \in C, \tag{5*}$$

and such that

$$y_* = Ax_* \in D \text{ solves } g(y_*)(y - y_*) \not\prec^V 0, \text{ for all } y \in D. \tag{6*}$$

We have the following special cases of nonlinear split ordered variational inequality problems:

1. When looked at the problems (5) and (6) separately, the problem (5) is the classical ordered variational inequality problem OVI($f, C, U$) on vector spaces. :

2. When we consider a special case of NSOVI($f, C, U, A, g, D, V$) such as: $X = Y$, $U = V$, $C = D$, $f = g$ and $A = I$, that is the identity mapping on $X$, then the nonlinear split ordered variational inequality problem NSOVI($f, C, U, A, g, D, V$) becomes the classical ordered variational inequality problem OVI($f, C, U$).

In next special case, we extend the linear split variational inequality problems defined in (1) and (2) to nonlinear split variational inequality problems.

3. Let $X$ and $Y$ be topological vector spaces with dual spaces $X^*$ and $Y^*$, respectively. Without any confusion caused, we use $\langle \cdot, \cdot \rangle$ to denote the pairings between $X^*$ and $X$ and between $Y^*$ and $Y$. For any given nonempty subsets $C \subseteq X$ and $D \subseteq Y$, given operators $f: C \to X^*$, $g: D \to Y^*$ and a nonlinear operator $A: X \to Y$, the nonlinear split variational inequality problem associated with $f, g, A, C, D$ (denoted by SVI($f, C, A, g, D$)) is formulated as follows:

$$\text{find } x_* \in C \text{ such that } \langle f(x_*), x - x_* \rangle \geq 0 \text{ for all } x \in C, \tag{7}$$

and such that

$$y_* = Ax_* \in D \text{ solves } \langle g(y_*), y - y_* \rangle \geq 0 \text{ for all } y \in D. \tag{8}$$

So from a different view, nonlinear split ordered variational inequality problems are natural extensions of classical ordered variational inequality problems in vector spaces and can be considered as natural extensions of nonlinear split variational inequality problems on vector spaces.

3.2 The main theorem

**Definition 3. 1.** *Let $(X, \succcurlyeq^X)$ and $(U, \succcurlyeq^U)$ be partially ordered vector spaces. Let $C$ be a nonempty subset of $X$. A mapping $f: C \to L(X, U)$ is said to be order-positive whenever, for any $x, y \in C$, if $x \preccurlyeq^X y$, then*

$$f(x)(z) \succcurlyeq^U 0 \quad \text{implies} \quad f(y)(z) \succcurlyeq^U 0, \text{ for } z \in X.$$

*In particular, if $(U, \succcurlyeq^U) = (R, \geq)$, a mapping $f: C \to X^*$ is order-positive whenever, for any $x, y \in C$, if $x \preccurlyeq^X y$, then*

$$\langle f(x), z \rangle \geq 0 \quad \text{implies} \quad \langle f(y), z \rangle \geq 0, \text{ for } z \in X.$$

Let the spaces $X, Y, (U, \succcurlyeq^U)$ and $(V, \succcurlyeq^V)$ and the mappings $f, g,$ and $A$ be given in Definition 3, we write

$$\pi(_fA_g)(x) = \{z \in C: f(x)(z) = \wedge\{f(x)(t): t \in C\}$$

$$\text{and } g(Ax)(Az) = \wedge\{g(Ax)(s): s \in D\}\}, \text{ for } x \in C. \tag{9}$$

Here without confusion, $\wedge\{f(x)(t): t \in C\}$ denotes the $\succcurlyeq^U$-smallest element of $\{f(x)(t): t \in C\} \subseteq (U, \succcurlyeq^U)$ and $\wedge\{g(Ax)(s): s \in D\}$ denotes the $\succcurlyeq^V$-smallest element of $\wedge\{g(Ax)(s): s \in D\} \subseteq (V, \succcurlyeq^V)$.

Now we prove the main theorem of this paper.

**Theorem 3.2.** *Let $(X, \succcurlyeq^X), (Y, \succcurlyeq^Y), (U, \succcurlyeq^U)$ and $(V, \succcurlyeq^V)$ be partially ordered vector spaces. Let $C$ be an $\succcurlyeq^X$-chain-complete subset of $X$ and $D$ a nonempty subset of $Y$. Let $f: C \to L(X, U), g: D \to L(Y, V)$ and $A: X \to Y$ be nonlinear mappings. If $f, g$ and $A$ satisfy the following conditions*:

V1. *$f$ and $g$ both are order-positive*;
V2. *For every $x \in C, \pi(_fA_g)(x)$ is a nonempty universally inductive subset of $C$*;
V3. *There are elements $x' \in C$ and $u' \in \pi(_fA_g)(x')$ satisfying $x' \preccurlyeq^X u'$*,
V4. *The operator $A$ is order-increasing*,

*then* NSOVI(*f, C, U, A, g, D, V*) *has a solution. Moreover*

   (i)   $(S(f, A, g), \geqslant^X)$ *is a nonempty inductive poset*;
   (ii)  $(S(f, A, g) \cap [x'), \geqslant^X)$ *is a nonempty inductive poset*.

*Consequently, we have*

   (i)′ NSOVI(*f, C, U, A, g, D, V*) *has an $\geqslant^X$-maximal solution*;
   (ii)′ NSOVI(*f, C, U, A, g, D, V*) *has an $\geqslant^X$-maximal solution $x^*$ with $x^* \geqslant^X x'$*.

*Proof.* The proof of this theorem is based on the fixed point theorem A applied to chain-complete subsets of partially ordered vector spaces.

Define a set-valued mapping $F$ on $C$ by

$$F(x) = \pi(_fA_g)(x)$$
$$= \{z \in C: f(x)(z) = \wedge\{f(x)(t): t \in C\} \text{ and } g(Ax)(Az) = \wedge\{g(Ax)(s): s \in D\}\}, \text{ for } x \in C.$$

From condition V2 in this theorem, $F: C \to 2^C \setminus \{\varnothing\}$ is a well-defined set-valued mapping. Next we show that $F$ is order-increasing upward. To this end, it is sufficient to show that, for any $x_1, x_2 \in C$, $x_1 \preccurlyeq^X x_2$ implies $F(x_1) \subseteq F(x_2)$. For any given $z \in F(x_1)$, we have

$$f(x_1)(z) \preccurlyeq^U f(x_1)(t), \text{ for all } t \in C \tag{10}$$

and

$$g(Ax_1)(Az) \preccurlyeq^V g(Ax_1)(s), \text{ for all } s \in D. \tag{11}$$

Since $f(x_1) \in L(X, U)$, from the order-linearity of the partially ordered vector spaces, and from (10), it implies that $f(x_1)(t - z) \geqslant^U 0$, for all $t \in C$. From the order-positive property of $f$, we get $f(x_2)(t - z) \geqslant^U 0$, for all $t \in C$. It implies $f(x_2)(z) \preccurlyeq^U f(x_2)(t)$, for all $t \in C$. That is,

$$f(x_2)(z) = \wedge\{f(x_2)(t): t \in C\}. \tag{12}$$

Since $g(Ax_1) \in L(Y, V)$, from (11), it implies that

$$g(Ax_1)(s - Az) \geqslant^V 0, \text{ for all } s \in D. \tag{13}$$

By condition V4 in this theorem, the assumption $x_1 \preccurlyeq^X x_2$ implies $Ax_1 \preccurlyeq^Y Ax_2$. From the order-positive property of $g$ and by (13), we get

$$g(Ax_2)(s - Az) \geqslant^V 0, \text{ for all } s \in D.$$

It implies $g(Ax_2)(Az) \preccurlyeq^V g(Ax_2)(s)$, for all $s \in D$. That is,

$$g(Ax_2)(Az) = \wedge\{g(Ax_2)(s): s \in D\}. \tag{14}$$

By combining (12) and (14) and from the definition $F$, it follows that $z \in F(x_2)$. Hence, $F$ is order-increasing upward, and $F$ satisfies condition A1 in Theorem A.

From condition V2, for every $x \in C$, $F(x)$ is a nonempty universally inductive subset of $C$. So $F$ satisfies condition A2 in Theorem A.

For the element $x'$ given in condition V3, it yields $u' \in F(x')$ and it satisfies $x' \preccurlyeq^X u'$. So condition A3 in Theorem A holds for this mapping $F$.

Then from Theorem A, $F$ has a fixed point on $C$. Next we show the following

$$\mathcal{F}(F) = \mathcal{S}(f, A, g). \tag{15}$$

For any $x^* \in \mathcal{F}(F)$, we have $x^* \in F(x^*)$. It implies that

$$f(x^*)(x^*) = \wedge\{f(x^*)(t): t \in C\} \text{ and } g(Ax^*)(Ax^*) = \wedge\{g(Ax^*)(s): s \in D\}$$

Hence the following two order inequalities simultaneously hold

$$f(x_*)(x - x_*) \succcurlyeq^U 0, \text{ for all } x \in C, \tag{5}$$

and

$$g(Ax^*)(y - Ax^*) \succcurlyeq^V 0, \text{ for all } y \in D. \tag{6}$$

Let $y^* = Ax^*$. Then $x^*$ (or $(x^*, y^*)$) is a solution of the problem NSOVI($f, C, U, A, g, D, V$). Hence we obtain

$$\mathcal{F}(F) \subseteq \mathcal{S}(f, A, g).$$

On the other hand, we can show that every solution to the problem NSOVI($f, C, U, A, g, D, V$) is a fixed point of $F$ on $C$, that is, $\mathcal{S}(f, A, g) \subseteq \mathcal{F}(F)$. Hence (15) is proved.

Then rest of this theorem follows immediately from Theorem A. □

**Remarks on Theorem 3.2**:

1. The considered input spaces $X$ and $Y$ in Theorem 3.2 are vector spaces, which are the most general underlying spaces for defining variational inequality problems;

2. The output spaces $(U, \succcurlyeq^U)$ and $(V, \succcurlyeq^V)$ Theorem 3.2 are partially ordered vector spaces, that generalizes the range sets from linear order (totally ordered) sets to partially ordered sets  Hence the results of Theorem 3.2 includes split vector variational inequality problems as a special case.

3. Since the underlying spaces do not have topological structures, the considered mappings are not required any continuity, which just satisfy some order-increasing conditions.

3.3 Applications to nonlinear split variational inequalities on partially ordered topological vector spaces

In this subsection, we consider nonlinear split variational inequality problems on partially ordered topological vector spaces NSVI($f, C, A, g, D$) defined in (7) and (9), in which $X$ and $Y$ are topological vector spaces with dual spaces $X^*$ and $Y^*$, respectively. In this case, the operator $\pi(_fA_g): C \to 2^C$ defined in (9) becomes the following operator:

$$m(_fA_g)(x) = \{z \in C: \langle f(x), z\rangle = \min\{\langle f(x_*), t\rangle: t \in C\}$$
$$\text{and } \langle g(Ax), Az\rangle = \min\{\langle g(Ax), s\rangle: s \in D\}\}, \text{ for } x \in C.$$

If the operator $\pi(_fA_g)$ is replaced by the new operator $m(_fA_g)$ in the proof of Theorem 3.2, then we immediately obtain the following result about the existence of solutions to nonlinear split variational inequality problems on partially ordered topological vector spaces.

**Corollary 3.3**. *Let $(X, \succcurlyeq^X)$ and $(Y, \succcurlyeq^Y)$ be partially ordered topological vector spaces with dual spaces $X^*$ and $Y^*$, respectively. Let C be a nonempty compact subset in X and D a nonempty subset in Y. Let $f: C \to X^*$, $g: D \to Y^*$ and $A: X \to Y$ be nonlinear mappings.*

*If f, g and A satisfy the conditions* V1, V4 *given in Theorem 3.2, and*

    V2'. *For every $x \in C$, $m(_fA_g)(x)$ is a nonempty closed subset of C;*
    V3'. *There are elements $x' \in C$ and $u' \in m(_fA_g)(x')$ satisfying $x' \preccurlyeq^X u'$,*

*then NSVI(f, C, A, g, D) has a solution and the solution set satisfies* (i), (ii), (i)', *and* (ii)' *in Theorem 3.2.*

*Proof.* From the listed examples of chain-complete posets in section 2, the compactness of $C$ implies that $(C, \succcurlyeq^X)$ is chain-complete. The condition that $m(_fA_g)(x)$ is a nonempty closed subset of the compact subset $C$ implies that $m(_fA_g)(x)$ is also compact; and therefore $m(_fA_g)(x)$ is universally inductive. Hence condition V2' in this corollary implies that $f$, $g$ and $A$ satisfy condition V2 in Theorem 3.2. □

3.4 Applications to nonlinear split variational inequalities on partially ordered Banach spaces

The study on nonlinear split variational inequalities on general partially ordered Banach spaces can be included in Corollary 3.3 as special cases of partially ordered topological vector spaces. By (4) and consequently, from Corollary 3.3, we have

**Corollary 3.4**. *Let $(X, \succcurlyeq^X)$ be a partially ordered Banach space with dual space $X^*$ and $(Y, \succcurlyeq^Y)$ a partially ordered topological vector spaces with dual spaces $Y^*$. Let C be a nonempty weakly compact subset in X and D a nonempty subset in Y. Let $f: C \to X^*$, $g: D \to Y^*$ and $A: X \to Y$ be nonlinear mappings.*

*If f, g and A satisfy the conditions* V1, V4 *given in Theorem 3.2, and*

    V2''. *For every $x \in C$, $m(_fA_g)(x)$ is a nonempty weakly compact subset of C;*
    V3'. *There are elements $x' \in C$ and $u' \in m(_fA_g)(x')$ satisfying $x' \preccurlyeq^X u'$,*

*then NSVI(f, C, A, g, D) has a solution and the solution set satisfies* (i), (ii), (i)', *and* (ii)' *in Theorem 3.2.*

*Proof.* Since the norm-topology on $X$ is natural with respect to the partial order $\succcurlyeq^X$ on $X$, so is the weak topology on $X$. It implies that $(C, \succcurlyeq^X)$ is chain-complete and, for every $x \in C$, $m(_fA_g)(x)$ is universally inductive. □

In particular, on partially ordered reflexive Banach spaces, similar to Corollary 3.4, we have the following result.

**Corollary 3.5**. *Let $(X, \succcurlyeq^X)$ be a partially ordered reflexive Banach space with dual space $X^*$ and $(Y, \succcurlyeq^Y)$ a partially ordered topological vector spaces with dual spaces $Y^*$. Let $C$ be a bounded, closed and convex subsets of $X$ and $D$ a nonempty subset of $Y$. Let $f: C \to X^*$, $g: D \to Y^*$ and $A: X \to Y$ be nonlinear mappings.*

*If $f$, $g$ and $A$ satisfy conditions* V1, V4 *given in Theorem 3.2, and* V2″, V3′ *in Corollary 3.4, then* NSVI$(f, C, A, g, D)$ *has a solution and the solution set satisfies* (i), (ii), (i)′, *and* (ii)′ *in Theorem 3.2.*

*Proof.* Since $C$ is a bounded, closed and convex subset of reflexive Banach space $X$, $C$ is weakly compact. Then $(C, \succcurlyeq^X)$ is chain-complete. □

For applications of Corollary 3.5, it is worth to note that condition V2″ in Corollary 3.5 can be replaced by:

V2‴. For every $x \in C$, $m(_fA_g)(x)$ is a nonempty closed and convex weakly subset of $C$.

## 4. (Linear) split variational inequalities on partially ordered Banach spaces

When looking at the condition V2 in Theorem 3.2 regarding to the mapping $\pi(_fA_g)$, condition V2′ in Corollary 3.3 and condition V2″ in Corollary 3.4 about the mapping $m(_fA_g)$, one may consider these conditions to be too strong and to be indirectly assumed for the mappings $f$, $g$, and $A$. In this section, we investigate some direct conditions about the mappings $f$, $g$, and $A$ for satisfaction of condition V2″ on partially ordered Banach spaces.

**Lemma 4.1**. *Let $(X, \succcurlyeq^X)$ and $(Y, \succcurlyeq^Y)$ be partially ordered topological vector spaces with dual spaces $X^*$ and $Y^*$, respectively. Let $C$ be a weakly compact subset of $X$ and $D$ a nonempty subset of $Y$. Let $f: C \to X^*$, $g: D \to Y^*$ be mappings and $A \in L(X, Y)$. Suppose that $f$, $g$ and $A$ satisfy the following conditions*

V5. $AC = D$;
V6. $A$ *is nonnegative reserving for $f$ and $g$, that is, for any $x, t \in C$,*

$$\langle f(x), t \rangle \geq 0 \quad \text{implies} \quad \langle g(Ax), At \rangle \geq 0.$$

*Then, for every $x \in C$,*

$$m(_fA_g)(x) = \{z \in C: \langle f(x), z \rangle = \min\{\langle f(x_*), t \rangle: t \in C\} \tag{16}$$

*and*

$$m(_fA_g)(x) \text{ is a nonempty weakly compact subset of } C.$$

*Proof.* For every $x \in C$, since $f(x) \in X^*$ and $C$ is a weakly compact subset of $X$, we have

$$\{z \in C: \langle f(x), z \rangle = \min\{\langle f(x_*), t \rangle: t \in C\} \neq \varnothing.$$

For any $z \in C$ satisfying $\langle f(x), z \rangle = \min\{\langle f(x_*), t \rangle: t \in C\}$, it follows

$$\langle f(x), t - z \rangle \geqslant 0, \text{ for all } t \in C.$$

That is,
$$\langle f(x), z \rangle = \min\{\langle f(x_*), t \rangle : t \in C\}. \tag{17}$$

From condition V6 and $A \in L(X, Y)$ in this lemma and by (17), it implies
$$\langle g(Ax), At - Az \rangle \geqslant 0, \text{ for all } t \in C. \tag{18}$$

From condition V5: $AC = D$ and by (18), we get
$$\langle g(Ax), s - Az \rangle \geqslant 0, \text{ for all } s \in D.$$

That is,
$$\langle g(Ax), Az \rangle = \min\{\langle g(Ax), s \rangle : s \in D\}\}. \tag{19}$$

Combining (17) and (19), it implies $z \in m(_fA_g)(x)$ and hence $m(_fA_g)(x) \neq \emptyset$. We obtain
$$\{z \in C : \langle f(x), z \rangle = \min\{\langle f(x_*), t \rangle : t \in C\} \subseteq m(_fA_g)(x).$$

It is clear that
$$m(_fA_g)(x) \subseteq \{z \in C : \langle f(x), z \rangle = \min\{\langle f(x_*), t \rangle : t \in C\}.$$

Hence
$$m(_fA_g)(x) = \{z \in C : \langle f(x), z \rangle = \min\{\langle f(x_*), t \rangle : t \in C\}, \text{ for } x \in C.$$

So (16) is proved. Next we prove that $m(_fA_g)(x)$ is weakly closed in $C$. To this end, from (16), we only need to prove that the set $\{z \in C : \langle f(x), z \rangle = \min\{\langle f(x_*), t \rangle : t \in C\}$ is weakly closed in $C$. Take a weakly convergent sequence $\{z_n\} \subseteq \{z \in C : \langle f(x), z \rangle = \min\{\langle f(x_*), t \rangle : t \in C\} = m(_fA_g)(x) \subseteq C$. Since $C$ is weakly compact, then there is $z_0 \in C$ such that
$$z_n \to z_0, \text{ weakly as } n \to \infty. \tag{20}$$

From $f(x) \in X^*$ and (20), it implies that,
$$\langle f(x), z_n \rangle \to \langle f(x), z_0 \rangle, \text{ as } n \to \infty. \tag{21}$$

By $\{z_n\} \subseteq \{z \in C : \langle f(x), z \rangle = \min\{\langle f(x_*), t \rangle : t \in C\}$ and from (21), we have
$$\langle f(x), z_0 \rangle = \min\{\langle f(x_*), t \rangle : t \in C\}.$$

Hence $z_0 \in \{z \in C : \langle f(x), z \rangle = \min\{\langle f(x_*), t \rangle : t \in C\}$. It follows that $m(_fA_g)(x)$ is weakly closed in this weakly compact set $C$. So $m(_fA_g)(x)$ is weakly compact. □

**Theorem 4.2.** *Let $(X, \geqslant^X)$ be a partially ordered Banach spaces with dual spaces $X^*$ and $(Y, \geqslant^Y)$ a partially ordered topological vector spaces with dual spaces $Y^*$. Let $C$ be a nonempty weakly compact subset in $X$ and $D$ a nonempty subset in $Y$. Let $f: C \to X^*$, $g: D \to Y^*$ be mappings and $A \in L(X, Y)$.*

*If $f$, $g$ and $A$ satisfy conditions V1, V4 given in Theorem 3.2, V3′ in Corollary 3.4 and V5, V6 in Lemma 3.4, then SVI($f$, $C$, $A$, $g$, $D$) has a solution and the solution set satisfies* (i), (ii), (i)′, *and* (ii)′ *in Theorem 3.2.*

*Proof.* Let $\omega$ be the weak topology on $X$. By remark (28), $(X, \omega, \geqslant^X)$ is a partially ordered topological space and $C$ is $\omega$-compact, which implies that $(C, \geqslant^X)$ is $\geqslant^X$-chain-complete. Then from Lemma 4.1, conditions V5, V6 imply that condition V2″ in Corollary 3.4 is satisfied. Then this corollary follows from Corollary 3.4 immediately. □

As a consequence of Theorem 4.2, we have

**Corollary 4.3**. *Let $(X, \geqslant^X)$ be a partially ordered reflexive Banach spaces with dual spaces $X^*$ and $(Y, \geqslant^Y)$ a partially ordered topological space with dual space $Y^*$. Let $C$ be a bounded, closed and convex subsets of $X$ and $D$ a nonempty subset of $Y$. Let $f: C \to X^*$, $g: D \to Y^*$ be mappings and $A \in L(X, Y)$.*

*If $f$, $g$ and $A$ satisfy conditions V1, V4 given in Theorem 3.2, V3′ in Corollary 3.4 and V5, V6 in Lemma 3.4, then $\text{SVI}(f, C, A, g, D)$ has a solution and the solution set satisfies* (i), (ii), (i)′, *and* (ii)′ *in Theorem 3.2*.

**References**

[1]   A. Auslender, "Optimization. Methodes Numeriques". Masson, Paris, 1976.

[2]   C. L. Byrne, A unified treatment of some iterative algorithms in signal processing and image reconstruction, *Inverse Problems*, **20** (2004), 103–120.

[3]   C. L. Byrne, Iterative oblique projection onto convex sets and the split feasibility problem, *Inverse Problems*, 18 (2002), 441–453.

[4]   Y. Censor, A. Gibali and S. Reich, Algorithms for split variational inequality problem, *Numerical Algorithms*, 59 (2012), 301–323.

[5]   Y. Censor, T. Elfving, A multi-projection algorithm using Bregman projections in a product space, *Numer. Algorithms*, 8 (1994), 221-239.

[6]   Y. Censor, T. Bortfeld, B. Martin and A. Tromov, A unified approach for inversion problems in intensity-modulated radiation therapy, *Physics in Medicine and Biology*, 51 (2006), 2353–2365.

[7]   Y. Censor and A. Segal, On string-averaging for sparse problems and on the split common Fixed point problem, *Contemporary Mathematics*, 513 (2010), 125–142.

[8]   S. S. Chang and Ravi P Agarwal, Strong convergence theorems of general split equality problems for quasi-nonexpansive mappings, *Journal of Inequalities and Applications* (2014), 344–367.

[9]   S. S. Chang, Lin Wang, Yong Kun Tang and Gang Wang, Moudafis open question and simultaneous iterative algorithm for general split equality variational inclusion problems and general split equality optimization problems, *Fixed Point Theory and Applications*, (2014), 301–215.

[10] H. Che and M. Li, A simultaneous iterative method for split equality problems of two finite families of strictly psuedononspreading mappings without prior knowledge of operator norms, *Fixed Point Theory Appl*. 2015, Article ID 1 (2015).

[11] K. Fan, A generalization of Tychonoff's fixed point theorem, *Math. Ann*. **142** (1961), 305–310.

[12] J. L. Li, Several extensions of the Abian–Brown fixed point theorem and their applications to extended and generalized Nash equilibria on chain-complete posets, *J. Math. Anal. Appl.*, (2014), 409:1084–1092.

[13] J. L. Li, Inductive properties of fixed point sets of mappings on posets and on partially ordered topological spaces, *Fixed Point Theory and Applications*, (2015) 2015:211 DOI 10.1186/s13663-015-0461-8.

[14] J. L. Li, Existence of Solutions to Split Variational Inequality Problems and Split Minimization Problems in Banach Spaces. It is available at:http://arxiv.org/abs/1511.07068

[15] J. L. Li and H. K. Xu, Solvability of Split Vector Variational Inequality Problems in Banach Spaces, to appear in JNCA.

[16] I. V., Konnov and J. C. Yao, On the generalized vector variational inequality problem, J. of Math. Anal. And Appl., **206** (1997), 42–58.

[17] A. Moudafis, Split monotone variational inclusions, *J. Optim. Theory Appl*. **150** (2011), 275–283.

[18] S. Park, Recent applications of Fan-KKM theorem, Lecture notes in Math. Anal. Institute, **1841** (2013), 58–68.

[19] W. Takahashi, The split feasibility problem in Banach spaces, *Journal of Nonlinear and Convex Analysis* **15** (2014), 1349–1355

[20] W. Takahashi, "Iterative methods for split common fixed point problems in Banach spaces", presented at the 11[th] International Conference on Fixed Point Theory and Its Applications, Istanbul, Turkey, July, 2015.

[21] Y. Xiao, Y. Censor, D. Michalski and J.M. Galvin, The least-intensity feasible solution for aperture-based inverse planning in radiation therapy, *Annals of Operations Research*, 119 (2003), 183–203.

[22] H. K. Xu, A variable Krasnosel'skii-Mann algorithm and the multiple sets split feasibility problem, *Inverse Problem*, 22 (2006), 2021-2034.